\documentclass{article}

\usepackage{amssymb}
\usepackage[english,francais]{babel}
\usepackage{amsthm}

\newcommand{\Q}{\mathbb{Q}}

\newcommand{\la}{\lambda}
\newcommand{\unsurr}{\frac{1}{r}}
\newcommand{\unsur}{\frac{1}}

\newcommand{\V}{\mathcal{V}}
\newcommand{\IH}{\mathcal{H}}

\newcommand{\lb}{\lbrace}
\newcommand{\rb}{\rbrace}
\newcommand{\noin}{\noindent}
\newcommand{\W}{\mathcal{W}}
\newcommand{\cil}{\frac{n(n-3)}{2}}
\newcommand{\dbw}{\frac{(n-1)(n-2)}{2}}
\newcommand{\chl}{\frac{n(n-1)}{2}}

\newcommand{\jca}{\unsur{r^{n-3}}}
\newcommand{\bmw}{\mathcal{B}}
\newcommand{\iwha}{\mathcal{H}}
\newcommand{\symn}{\mathcal{S}}

\title{Irreducibility of the Lawrence-Krammer representation of the BMW algebra of type
$A_{n-1}$}
\begin{document}
\selectlanguage{english}
\author{Claire I. Levaillant\\\\\textit{Caltech, MC 253-37, Pasadena, CA
91125, cl@caltech.edu}}
\date{}
\maketitle

\begin{abstract}
\selectlanguage{english} \noindent The Lawrence-Krammer
representation introduced by Lawrence and Krammer in order to show
the linearity of the braid group is generically irreducible. We show
this fact and show further that for some values of its two
parameters, when these are specialized to complex numbers, the
representation becomes reducible. We describe what these values are
and give a complete description of the dimensions of the invariant
subspaces when the
representation is reducible. \\\\

\selectlanguage{francais} \begin{center}{\bf R\'esum\'e}\end{center}
\noindent La repr\'esentation de Lawrence-Krammer, introduite par
Lawrence et Krammer pour montrer la lin\'earit\'e du groupe de
tresses, est g\'en\'eriquement irr\'eductible. On montre ce fait et
on montre \'egalement que lorsque les deux param\`etres de la
repr\'esentation prennent certaines valeurs complexes, la
repr\'esentation devient r\'eductible. On donne ici toutes les
valeurs des param\`etres pour lesquelles la repr\'esentation est
r\'eductible, ainsi que les dimensions des sous-espaces stables.

\end{abstract}

\selectlanguage{francais}
\section*{Version fran\c{c}aise abr\'eg\'ee}
\noindent On consid\`ere, suivant la d\'efinition de \cite{CGW},
l'alg\`ebre de Birman-Murakami-Wenzl $\mathcal{B}$ de type $A_{n-1}$
avec param\`etres complexes non nuls $l$ et $m$ sur le corps
$\Q(l,r)$, o\`u les param\`etres $r$ et $m$ sont reli\'es par
$m=\unsurr-r$. Cette alg\`ebre a pour g\'en\'erateurs
$g_1,\dots,g_{n-1}$, qui satisfont aux relations de tresses et elle
contient aussi d'autres \'el\'ements $e_1,\dots,e_{n-1}$, qui sont
d\'efinis par les relations:
$$ e_i=\frac{l}{m}(g_i^2+m\,g_i-1)\qquad i=1,\dots,n-1$$ On travaille
sur une repr\'esentation de degr\'e $\chl$ de cette alg\`ebre dans
l'espace de Lawrence-Krammer. Comme repr\'esentation du groupe de
tresses \`a $n$ brins, cette repr\'esentation est \'equivalente, \`a
r\'e\'echelonnage pr\`es des g\'en\'erateurs, \`a la
repr\'esentation du groupe d'Artin de type $A_{n-1}$ bas\'ee sur les
param\`etres $r$ et $t$, d\'ecrite dans \cite{CW}. Cette derni\`ere
repr\'esentation est elle-m\^eme \'equivalente \`a la
repr\'esentation de Lawrence-Krammer bas\'ee sur les param\`etres
$q$ et $t$ et utilis\'ee par Krammer et ind\'ependamment par Bigelow
pour montrer la lin\'earit\'e du groupe de tresses \`a $n$ brins
(Voir \cite{KR} et \cite{BIG}). Le $r$ de cette note est le
$\unsurr$ de \cite{CGW} et les deux param\`etres $r$ et $l$ de
\cite{CGW} sont reli\'es aux param\`etres $t$ et $r$ de \cite{CW}
par $lt=\unsur{r^3}$. Les param\`etres $t$ et $r$ de \cite{CW} sont
eux-m\^eme reli\'es aux param\`etres $t$ et $q$ de \cite{KR} par
$q=r^2$ (cf \cite{CW}, exemple $3.8$). L'irr\'eductibilit\'e de la
repr\'esentation de Lawrence-Krammer pour des valeurs g\'en\'eriques
de $q$ et $t$ est contenue dans les travaux de Zinno (voir
\cite{Z}), de Cohen-Gijsbers-Wales (voir \cite{CGW}) et de Marin
(voir \cite{MA2}). D'autres travaux \'etudient aussi cette
repr\'esentation et traitent de sa r\'eductibilit\'e. En
particulier,
certaines conditions suffisantes d'irr\'eductibilit\'e sont
implicitement dans la th\`ese de Marin (voir \cite{MA1}) et dans
\cite{B}, Bigelow \'etudie le cas de r\'eductibilit\'e
$t=\unsur{q}$. Le th\'eor\`eme suivant donne toutes les valeurs
complexes des param\`etres pour lesquelles la repr\'esentation est
r\'eductible et
d\'ecrit les dimensions des sous-espaces stables irr\'eductibles.\\

\newtheorem*{theo}{Th\'eor\`eme}
\begin{theo}
On se donne trois param\`etres complexes non nuls $l$, $m$ et $r$,
o\`u $m=\unsurr-r$ et un entier $n\geq 3$ et on suppose que
l'alg\`ebre de Iwahori-Hecke du groupe sym\'etrique $\mathcal{S}_n$
avec param\`etre $r^2$ sur le corps $\Q(l,r)$ est semisimple. \\\\
$\bullet$ La repr\'esentation de Lawrence-Krammer de l'alg\`ebre BMW
$\mathcal{B}(A_{n-1})$ avec param\`etres $l$ et $m$ sur le corps
$\Q(l,r)$ est irr\'eductible, sauf quand:
\begin{list}{\texttt{*}}{}
\item $l=r$ et $n\geq 4$: il existe alors un unique sous-espace stable irr\'eductible de dimension
$\cil$ dans l'espace de Lawrence-Krammer.
\item $l=-r^3$: il existe alors un sous-espace stable irr\'eductible de dimension
$\dbw$ dans l'espace de Lawrence-Krammer. De plus, celui-ci est
unique sauf pour $n=3$ et $r^6=-1$.
\item $l=\unsur{r^{2n-3}}$: il existe alors un sous-espace stable
unidimensionel dans l'espace de Lawrence-Krammer. De plus, celui-ci
est unique, sauf pour $n=3$ et $r^6=-1$.
\item $l=\unsur{r^{n-3}}$ ou $l=-\unsur{r^{n-3}}$: il existe alors
un unique sous-espace stable irr\'eductible de dimension $(n-1)$
dans l'espace de Lawrence-Krammer.
\end{list}\vspace{0.03in}
$\bullet$ Quand $n=3$ et $l=-r^3=\unsur{r^3}$, il existe exactement
deux sous-espaces stables unidimensionnels. \\ $\bullet$ Les
sous-espaces stables irr\'eductibles mentionn\'es sont les seuls
sous-espaces stables irr\'eductibles pouvant appara\^itre dans
l'espace de Lawrence-Krammer. \\ $\bullet$ Pour ces valeurs des
param\`etres, la repr\'esentation est r\'eductible et
ind\'ecomposable, donc
l'alg\`ebre BMW n'est pas semisimple.\\

On se donne deux param\`etres complexes non nuls $t$ et $q$, o\`u
$q$ n'est pas une racine $k$-i\`eme de l'unit\'e pour tout entier
$k$ tel que $1\leq k\leq n$. La repr\'esentation de Lawrence-Krammer
du groupe de tresses $B_n$ \`a $n$ brins, bas\'ee
sur les param\`etres $q$ et $t$ est r\'eductible si et seulement si\\
$$\left.\begin{array}{l}t\in\lb \unsur{q},-1,\unsur{q^n},\unsur{\sqrt{q}^n},-\unsur{\sqrt{q}^n}\rb
\;\;\;quand\;\;\; n\geq 4\\ t\in\lb
-1,\unsur{q^3},\unsur{\sqrt{q}^3},-\unsur{\sqrt{q}^3}\rb
\;\;\;\;\;\;\;\,dans\;\;\; le\;\;\; cas\;\;\; de\;\;\;
B_3\end{array}\right.$$ \end{theo}
\selectlanguage{english}
\section{The representation}

Let $l$ and $m$ be two nonzero complex parameters and let $r$ and
$-\unsurr$ be the two complex roots of the quadratics $X^2+mX-1$.
Let $n$ be an integer with $n\geq 3$. Throughout the note, we will
assume that the Iwahori-Hecke algebra $\IH_n$ of the Symmetric group
$\symn_n$ with parameter $r^2$ over the field $\Q(l,r)$ is
semisimple. This assumption is met exactly when $r^{2k}\neq 1$ for
every integer $k$ with $1\leq k\leq n$ (See \cite{M}). We consider
the BMW algebra $\bmw$ of type $A_{n-1}$ with parameters $l$ and $m$
over the field $\Q(l,r)$, as defined in \cite{CGW}. This algebra has
two sets of $n-1$ elements: the $g_i$'s, $i=1,\dots,n-1$, that
satisfy the braid relations and generate $\bmw$ and the $e_i$'s,
$i=1,\dots,n-1$, that are related to the $g_i$'s by
$$e_i=\frac{l}{m}(g_i^2+m\,g_i-1)$$ We define $\iwha$ to be the Iwahori-Hecke algebra
of the symmetric group $\symn_{n-2}$ over the field $\Q(l,r)$ with
generators $g_3,\dots,g_{n-1}$ that satisfy the braid relations and
the relation $g_i^2+m\,g_i=1$ for all $i$. The base field $\Q(l,r)$
is an $\iwha$-module for the action $g_i.1=r$, where
$i=3,\dots,n-1$. We define $\bmw_1$ to be the $\bmw$-module-$\iwha$:
$$\bmw e_1/<\bmw e_ie_1|i=3,\dots,n-1>$$
When $n=3$, $\bmw_1$ is simply $\bmw e_1$. We now obtain a
$\bmw$-module of dimension $\chl$ over $\Q(l,r)$ by considering the
tensor product:
$$\V^{(n)}=\bmw_1\otimes_{\iwha}\Q(l,r)$$
This $\bmw$-module is precisely the generically irreducible
representation that we study throughout this note. As a matter of
fact, each product of the algebra $e_ie_j$, for non-adjacent nodes
$i$ and $j$, acts trivially on $\V^{(n)}$. Then, by \cite{CGW}, this
left representation of $\bmw$
must be equivalent to the Lawrence-Krammer representation of the BMW
algebra of type $A_{n-1}$. As a representation of the braid group on
$n$ strands, it is equivalent to the representation based on the two
parameters $q$ and $t$, introduced by Krammer in \cite{KR} to show
that the braid group on $n$ strands is linear. The link between the
parameters $l$ and $r$ of this paper and the parameters $q$ and $t$
of Krammer's representation is given by $q=\unsur{r^2}$ and
$lt=r^3$.

\section{Reducibility of the Representation}

Our main Theorem is the following. Details of the proof appear in \cite{CIL}.\\
\newtheorem{theorem}{Theorem}
\begin{theorem} Let $l$, $m$ and $r$ be three nonzero complex parameters,
where $m$ and $r$ are related by $m=\unsurr-r$. Assume that the
Iwahori-Hecke algebra $\IH_n$ of the symmetric group $\symn_n$ with
parameter $r^2$ over the field $\Q(l,r)$ is semisimple, that is
assume that $r^{2k}\neq 1$ for every integer $k\in\lb
1,\dots,n\rb$.\\\\ When $n\geq 4$, the Lawrence-Krammer
representation of the $BMW$ algebra $\mathcal{B}(A_{n-1})$ with
parameters $l$ and $m$ over the field $\Q(l,r)$ is irreducible,
except when $l\in\lb
r,-r^3,\unsur{r^{2n-3}},\unsur{r^{n-3}},-\unsur{r^{n-3}}\rb$ when it
is reducible.\\
When $n=3$, the Lawrence-Krammer representation of the $BMW$ algebra
$\mathcal{B}(A_2)$ with parameters $l$ and $m$ over the field
$\Q(l,r)$ is irreducible except when $l\in\lb
-r^3,\unsur{r^3},1,-1\rb$ when it is reducible.\\
\end{theorem}
\noin\textbf{Proof.} (Sketch) We show that if there exists a proper
invariant subspace $\W$ of $\V^{(n)}$, then all the $e_i$'s and
$g_k$ conjugates of the $e_i$'s, $i=1,\dots,n-1$, must annihilate
$\W$. In particular, the action of $\bmw$ on $\W$ is an
Iwahori-Hecke algebra action. For the small values
$n\in\lb3,4,5,6\rb$, we choose a basis of $\V^{(n)}$ and compute the
matrix $M(n)$ of the left action of the BMW algebra element
$$\sum_{1\leq i\leq n-1}e_i + \sum_{1\leq i<j-1<
n}g_{j-1}^{-1}\dots g_{i+1}^{-1}e_ig_{i+1}\dots g_{j-1}$$ in this
basis. If $\W$ is nontrivial, the determinant of this matrix must be
zero. By using Maple, it yields the values for $l$ and $r$ that are
described in the Theorem. Conversely, for the values of $l$ and $r$
that annihilate this determinant, we find some nonzero vectors in
the kernel $K(n)$ of the matrix $M(n)$. We show that $K(n)$ is a
proper submodule of $\V^{(n)}$. It then gives the equivalence of
Theorem $2.1$ for these small values of $n$. When $n$ is large
enough, that is when $n\geq 7$, the irreducible representations of
$\iwha_n$ have degrees $1$, $n-1$, $\cil$, $\dbw$ or degrees greater
than $\dbw$, except in the case $n=8$, when they have degrees $1$,
$7$, $14$, $20$ or $21$ (see \cite{J} and \cite{M}). We show that
there exists a one-dimensional invariant subspace of $\V^{(n)}$ if
and only if $l=\unsur{r^{2n-3}}$ and that there exists an
irreducible $(n-1)$-dimensional invariant subspace of $\V^{(n)}$ if
and only if $l\in\lb\unsur{r^{n-3}},-\unsur{r^{n-3}}\rb$. We then
proceed by induction on $n$. First, let $n=7$ or $n\geq 9$ and
suppose Theorem $2.1$ holds for representations of
$\mathcal{B}(A_k)$ where $k\in\lb n-3,n-2\rb$. Suppose that there
exists a proper invariant subspace $\W$ of $\V^{(n)}$ of dimension
greater than or equal to $\cil$. Then the intersections
$\W\cap\V^{(n-1)}$ and $\W\cap\V^{(n-2)}$ must be nontrivial. By
induction, $l$ must then belong to $$\bigg\lb
r,-r^3,\unsur{r^{2n-5}},\unsur{r^{n-4}},-\unsur{r^{n-4}}\bigg\rb\cap\bigg\lb
r,-r^3,\unsur{r^{2n-7}},\unsur{r^{n-5}},-\unsur{r^{n-5}}\bigg\rb$$
We assumed that $r^{2k}\neq 1$ for every $k=1,\dots,n$. Then, by
inspection, we see that the only possibility is to have $l\in\lb
r,-r^3\rb$. Let's deal with the case $n=8$. The case $n=8$ is in
fact not different. Indeed, if $\W$ has dimension greater than or
equal to $14$, then again, the spaces $\W\cap\V^{(7)}$ and
$\W\cap\V^{(6)}$ are nontrivial. Thus, we have proven that if the
representation is reducible, it forces the values of the theorem for
$l$ and $r$. It remains to show that when $l=r$ or $l=-r^3$, the
representation is reducible. In each case, we verify that the
nonzero vector belonging to $K(5)\cap\V^{(4)}$ that we found with
Maple also belongs to all the $K(n)$'s for $n\geq 6$.
This shows the reducibility of the representation in both cases.\\
\newtheorem{remarque}{\it Remark}

\begin{remarque}
In \cite{CIL}, we also give a complete proof without using Maple for
the small cases $n\in\lb3,4,5,6\rb$.
\end{remarque}
\begin{remarque}
In \cite{W}, Hans Wenzl states that $\bmw(A_{n-1})$ is semisimple,
except possibly if $r$ is a root of unity or $l$ is some power of
$r$, where he also considers complex parameters. Here Theorem $2.1$
and the method that we use imply that for these specific values of
$l$ and $r$, the algebra is not semisimple as the representation is
then reducible and indecomposable.
\end{remarque}
\section{Dimensions of the invariant subspaces when the representation is reducible}

In this section we give a series of theorems on the dimensions of
the invariant subspaces when the representation is reducible. We
still assume that $\IH_n$ is semisimple. We have the following results:\\


\begin{theorem} Let $n$ be an integer with $n\geq 4$. There exists a
one-dimensional invariant subspace of $\V^{(n)}$ if and only if
$l=\unsur{r^{2n-3}}$. If so, it is unique.\\
\textit{(Case $n=3$)} There exists a one-dimensional invariant
subspace of $\V^{(3)}$ if and only if $l\in\lb -r^3,\unsur{r^3}\rb$.
If $r^6\neq -1$, it is unique. If $r^6=-1$, there exists exactly two
one-dimensional invariant subspaces of $\V^{(3)}$.\\
\end{theorem}

\begin{theorem}
Let $n$ be an integer with $n\geq 3$ and $n\neq 4$. There exists an
irreducible $(n-1)$-dimensional invariant subspace of $\V^{(n)}$ if
and only if
$l\in\lb\unsur{r^{n-3}},-\unsur{r^{n-3}}\rb$. If so, it is unique.\\
\textit{(Case $n=4$)} There exists a $3$-dimensional invariant
subspace of $\V^{(4)}$ if and only if $l\in\lb\unsurr,
-\unsurr,-r^3\rb$. If so, it is unique.\\ \end{theorem}

\begin{theorem}
Let $n$ be an integer with $n\geq 4$. There exists an irreducible
$\cil$-dimensional invariant subspace of $\V^{(n)}$ if and only if
$l=r$. If so, it is unique.\\
\end{theorem}

\begin{theorem} Let $n$ be an integer with $n\geq 5$. There exists an
irreducible $\dbw$-dimensional invariant subspace of $\V^{(n)}$ if
and only if $l=-r^3$. If so, it is unique.
\\
\end{theorem}

\begin{theorem} These are the only irreducible invariant subspaces that may
appear in $\V^{(n)}$.\\
\end{theorem}
\noin\textbf{Proof of Theorem 2} (Sketch) Theorem $2$ is easily
proven by hand. Indeed, if $v$ is a spanning vector of a
one-dimensional invariant subspace of $\V^{(n)}$, and if
$\la_1,\dots,\la_{n-1}$ are scalars such that $g_i.v=\la_i\,v$, then
since the $e_i$'s all annihilate $v$, we get $\la_i^2+m\,\la_i-1=0$
for each $i$, so that the $\la_i$'s must take the values $r$ or
$-\unsurr$. Moreover, by the braid relations
$g_ig_{i+1}g_i=g_{i+1}g_ig_{i+1}$ for every node $i\in\lb
1,\dots,n-2\rb$, and the fact that $(r^2)^2\neq 1$, we see that all
the $\la_i$'s must take the same value $r$ or $-\unsurr$. Except
when $n=3$, where both values are possible, each of which
respectively forcing $l=\unsur{r^3}$ and $l=-r^3$, in the case when
$n\geq 4$, the $\la_i$'s must all equal $r$. This forces the value
$\unsur{r^{2n-3}}$ for $l$ as well as a unique (up to multiplication
by a scalar) spanning vector.\\

\noin\textbf{Proof of Theorem 3} (Sketch) When $n\geq 4$ and $n\neq
6$ (resp $n=6$) there are exactly two (resp four) inequivalent
irreducible representations of $\IH_n$ of degree $(n-1)$. When
$n\geq 5$ and $n\neq 6$ (resp $n=6$), we show that only one of these
two (resp four) representations may occur inside $\V^{(n)}$.
Moreover, this representation occurs when $l\in\lb\jca,-\jca\rb$ and
leads in each case to a unique irreducible invariant subspace of
$\V^{(n)}$ of dimension $(n-1)$. When $n=4$, both may occur, the
first one leading to the values $l\in\lb\unsurr,-\unsurr\rb$ and the
conjugate one forcing the value $l=-r^3$. For each of the three
(distinct) values of $l$, we can show the uniqueness of the
$3$-dimensional irreducible invariant subspace of $\V^{(4)}$. When
$n=3$, the irreducible representation of $\IH_3$ of degree $2$ is
self-conjugate and occurs in $\V^{(3)}$ for $l\in\lb 1,-1\rb$.
Again, for each of these values, there is a unique
irreducible $2$-dimensional invariant subspace in $\V^{(3)}$.\\

\noin\textbf{Proof of Theorem 4} (Sketch) We deal with the case
$n=4$ directly by hand. When $n\geq 5$, we show that when $l=r$, the
rank of the matrix $M(n)$ is greater than or equal to $n$. This is
simply achieved by exhibiting an invertible submatrix of $M(n)$ of
size $n$. If $k(n)$ denotes the dimension of the kernel $K(n)$, this
then yields the inequality on the dimensions $k(n)\leq\cil$. Also,
we show that the $\bmw$-module $K(n)$ is irreducible. Let's recall
from above that any irreducible proper invariant subspace of
$\V^{(n)}$ is an irreducible $\IH_n$-module. When $n\geq 5$ and
$n\neq 8$, it thus has dimension $1$, $(n-1)$, $\cil$, $\dbw$ or
dimension greater than $\dbw$ by \cite{J} and \cite{M}. This fact on
the dimensions and Theorems $2$ and $3$ imply that the dimension
$k(n)$ of $K(n)$ is in fact greater than or equal to $\cil$. The
latter inequality also holds for $n=8$, but this case needs to be
dealt with separately. Details can be found in \cite{CIL},
$\S\,9.2$. Thus, when $l=r$, the $\bmw$-module $K(n)$ is an
irreducible invariant subspace of $\V^{(n)}$ of dimension $\cil$.
Since from our discussion at the beginning of part $2$, any proper
invariant subspace of $\V^{(n)}$ must be contained in $K(n)$, the
uniqueness is established.

Conversely, we show that the existence of an irreducible
$\cil$-dimensional invariant subspace of $\V^{(n)}$ implies that
$l=r$ and we take the following path. First, we examine the case
$n=5$ by hand. We then proceed by induction on $n$. Except when
$n=7$, there are exactly two inequivalent irreducible
representations of $\IH_n$ of degree $\cil$. By inspection for
$n=5$, only one of them may occur inside $\V^{(5)}$. This is part of
our induction hypothesis. The proof is then completed by using the
branching rule. The case $n=7$ is in fact not an exception and also
follows from the branching rule, after showing that the Specht
modules $S^{(3,3)}$ and $S^{(2,2,2)}$ cannot occur in $\V^{(6)}$.
For details, see \cite{CIL}, $\S\,8.3$.\\

\noin\textbf{Proof of Theorem 5} (Sketch) First, if there exists an
irreducible $\dbw$-dimensional invariant subspace $\W$ of
$\V^{(n)}$, for $n\geq 6$, its dimension is large enough to make the
intersections $\W\cap\V^{(n-1)}$ and $\W\cap\V^{(n-2)}$ nontrivial
and to thus force $l\in\lb r,-r^3\rb$. Since when $l=r$, we know
that $k(n)=\cil$, which implies in particular that any proper
invariant subspace of $\V$ must have dimension less than or equal to
$\cil=\dbw-1$, we conclude that it is impossible to have $l=r$.
Thus, $l=-r^3$. The case $n=5$ needs to be dealt with separately.
Details appear in \cite{CIL}, $\S\,9.2$.

Conversely, suppose $l=-r^3$ and let $n\geq 5$. There are two cases:
\begin{list}{\texttt{*}}{}
\item If $r^{2n}\neq -1$, then we can show that $K(n)$ is
irreducible (this fact also holds for $n=4$). Moreover, we can show
that the rank of the matrix $M(n)$ is greater than or equal to
$(n-1)$, by exhibiting an invertible submatrix of $M(n)$ of size
$(n-1)$. Then $k(n)\leq\dbw$. By the previous studies it forces
$k(n)=\dbw$, the case $n=8$ being dealt with separately (see
\cite{CIL}, $\S\,9.3$). Hence there exists an irreducible $\dbw$
invariant subspace of $\V^{(n)}$.

\item If $r^{2n}=-1$, then $K(n)$ is no longer irreducible, but $r^{2(n-1)}\neq -1$, so that
$K(n-1)$ is irreducible. Still by the previous case when $n\geq 6$
and by Theorem $3$ when $n=5$, we have
$k(n-1)=\frac{(n-2)(n-3)}{2}$. By the discussion at the end of part
$2$, we know that $K(n)\cap\V^{(n-1)}$ is nontrivial. It follows
that
$$K(n)\cap\V^{(n-1)}=K(n-1),$$ by irreducibility of $K(n-1)$.
This yields on the dimensions: $k(n)\leq k(n-1)+(n-1)$. Also, we
have $K(n-1)\subset K(n)$. Further, we can show that: $k(n)\geq
k(n-1)+(n-2)$. With $k(n-1)=\frac{(n-2)(n-3)}{2}$, both inequalities
now yield $k(n)\in\big\lb\dbw,\dbw+1\big\rb$. We then get
$k(n)=1+\dbw$, as otherwise, the existence of a one-dimensional
invariant subspace would force the existence of an irreducible
$\cil$-dimensional invariant subspace by semisimplicity of $\IH_n$.
By Theorem $4$, $l$ would then be r, which contradicts our
assumption $l=-r^3$. From there, it appears that there exists an
irreducible $\dbw$-dimensional invariant subspace of $\V^{(n)}$. We can show that it is unique.\\
\end{list}

\begin{remarque} Joint proof of Theorems $4$ and $5$ is given
in \cite{CIL}, where we did not use the branching rule.
\end{remarque}

\end{document}